\documentclass{ifacconf}
\usepackage{xcolor}

\providecommand{\com[2]}{\begin{tt}[#1: #2]\end{tt}}

\usepackage{amsmath, amsfonts}
\newcommand{\xd}{x^{[d]}}
\newcommand{\dxd}{(x+\Delta x)^{[d]}}

\newcommand{\Axd}{(\mathbf{A_j}x)^{[d]}}
\newcommand{\Axdx}{(\mathbf{A_j}(x+\Delta x))^{[d]}}
\newcommand{\xk}{x^{\otimes d}}
\newcommand{\Ad}{\mathbf{A_j}^{[d]}}
\usepackage{graphicx}      
\usepackage{natbib}        
\begin{document}
\begin{frontmatter}
\title{Data-driven stability analysis of switched linear systems with Sum of Squares guarantees
} 
\thanks[footnoteinfo]{R. Jungers is a FNRS honorary Research Associate. This project has received funding from the European Research Council (ERC) under the European Union's Horizon 2020 research and innovation programme under grant agreement No 864017 - L2C. RJ is also supported by the Walloon Region, the Innoviris Foundation, and the FNRS (Chist-Era Druid-net).}

\author[First]{Anne Rubbens} 
  \author[First]{Zheming Wang} 
  \author[First]{Raphael M. Jungers}

  \address[First]{The ICTEAM Institute, UCLouvain, Louvain-la-Neuve,1348, Belgium
  (email: anne.rubbens@student.uclouvain.be, zheming.wang@uclouvain.be, raphael.jungers@uclouvain.be).}

\begin{abstract}                
We present a new data-driven method to provide probabilistic stability guarantees for black-box switched linear systems. By sampling a finite number of observations of trajectories, we construct approximate Lyapunov functions and deduce the stability of the underlying system with a user-defined confidence. The number of observations required to attain this confidence level on the guarantee is explicitly characterized.\\
Our contribution is twofold: first, we propose a novel approach for common quadratic Lyapunov functions, relying on sensitivity analysis of a quasi-convex optimization program. By doing so, we improve a recently proposed bound. Then, we show that our new approach allows for extension of the method to Sum of Squares Lyapunov functions, providing further improvement for the technique. We demonstrate these improvements on a numerical example. 
\end{abstract}

\begin{keyword}
Stability analysis of switched linear systems, Quasi-convex optimization, Data-driven analysis, Scenario approach, Sum of Squares techniques. 
\end{keyword}

\end{frontmatter}
\section{Introduction}
Hybrid systems, whose dynamics include continuous and discrete parts, allow to model a wide range of complex systems, including cyber-physical systems which have been a subject of growing interest recently (\cite{Lee:15}, \cite{Tab:14}). An important family of hybrid systems are \emph{switched linear systems}, consisting of several linear continuous state models, the \emph{modes}, together with a discrete rule deciding the switching among them. Formally, let $M=\{1,2,...,m\}$, and $\mathcal{M}=\{A_i, \ i \in M\}$ be a finite set of $m$ matrices in $\mathbb{R}^{n \times n}$, the modes. We consider systems of the form
\begin{equation}
    x_{k+1}=A_{\tau(k)}x_k
    \label{system}
\end{equation}
where $x_k \in \mathbb{R}^n$ is the state, $k \in \mathbb{N} $ is the time index, and $\tau \in M^{\mathbb{N}}$ is the \emph{switching sequence} taking values in $M$. Switched systems can for instance be used to model linear systems subject to time-dependent uncertainty or multi-controller switching schemes (\cite{Lib:03}, \cite{Sho:07}). 
\\
A system of the form (\ref{system}) is \emph{uniformly asymptotically stable} if it possesses the following property:
$$
    \forall \tau \in M^{\mathbb{N}}, \ \forall x_0 \in \mathbb{R}^n, \ ||x_k||\underset{k\longrightarrow \infty}{\longrightarrow}0.
$$
This property is characterized by the \emph{joint spectral radius} (JSR) of $\mathcal{M}$  (\cite{Jun:09}):
$$\rho (\mathcal{M})=\underset{k\longrightarrow\infty}{\lim \ } \underset{i_1,...,i_k}{\max \ } \{||A_{i_1}...A_{i_k}||^{\frac{1}{k}}: A_{i_j} \in \mathcal{M}\}.$$
Indeed, given a finite set of matrices $\mathcal{M}$, the corresponding switched dynamical system is uniformly asymptotically stable if and
only if $\rho(\mathcal{M})<1$ (\cite{Jun:09}, Cor. 1.1).
\\
Because switched linear systems do not inherit the properties of their modes, their stability analysis is difficult: deciding whether their JSR is smaller than 1 is NP-hard (\cite{Blo:99}), and many tools have therefore been designed in the past decades for their analysis. Most of these tools are model-based, build on notions from Lyapunov and invariant set theory, and provide deterministic
guarantees when the number of modes is finite (\cite{Lin:09},  \cite{Blan:08}, \cite{Lib:04}).
\\
However, in many practical applications, a closed-form model may not be available (think for example of self-driving cars, complex robotic systems, or smart grids applications). Hence, stability analysis based on information obtained
via a sufficient number of simulations (or observations of the unknown system) has received increasing attention lately (\cite{Koz:16}, \cite{Blan:17}, \cite{Ken:19}, \cite{Sma:20}, \cite{Wang:20}). A first approach consists in using model-based stability analysis techniques after recovering the
model from the simulations. However, since the identification problem for switched systems is NP-hard (\cite{Lau:16}), we bypass it in our approach.
We refer to our framework as \emph{data-driven analysis} of \emph{black-box} dynamical systems. It has been increasingly understood that one can obtain effective control methods, such as Lyapunov functions or invariant sets, despite the complexity of the underlying system and the lack of a closed- form description of the model.
\\
In \cite{Ken:19}, a probabilistic stability guarantee for black-box switched linear systems is provided, based on \emph{Common Quadratic Lyapunov Functions} (CQLFs) decreasing on a finite number $N$ of observations of their trajectories. The method proceeds by deriving a guarantee on the measure of the subset of the state space on which the CQLF is not decreasing using scenario optimization theory (\cite{Cam:18}, \cite{Cal:10}), and then transforming this guarantee into a global stability certificate with a certain confidence level, using geometric properties of switched systems.
\\
Relying on CQLFs is conservative for switched linear systems, since some stable switched systems do not admit a CQLF (\cite{Lin:09}). For this reason, in this paper, we derive a broader class of probabilistic stability guarantees for black-box switched linear systems. We generalize the data-driven analysis technique to a more complex class of Lyapunov functions decreasing along the observed trajectories, namely \emph{Sum of Squares (SOS) Lyapunov functions} (\cite{Par:08}, \cite{Jun:09}, Section 2.3.7). 
\\ 
Our main contribution consists in the derivation of a probabilistic upper bound on the JSR of a black-box system, with a user-defined confidence level, via a SOS Lyapunov function of arbitrary degree $2d$. The trade-off between the number of observations $N$ and the tightness of the upper bound is explicitly quantified. Because the geometric tools used in \cite{Ken:19} do not generalize in a straightforward way to SOS Lyapunov functions of higher degrees, we develop another approach, based on \emph{sensitivity analysis} of a quasi-convex optimization problem (\cite{Ber:16}). Hence, our new approach does not rely on classical chance-constrained theorems usually underlying the scenario approach. Additionally, we thereby also obtain a novel upper bound in the quadratic case, that can be compared with the one obtained in \cite{Ken:19}. In the white-box setting, approximations of the JSR using SOS Lyapunov functions are tighter (\cite{Par:08}), and we will observe that this generalizes to the black-box setting.
\\
The outline of this paper is as follows. This section ends with the notation, and Section 2 contains preliminary results about the approximation of the JSR using CQLFs as well as the introduction of the problem studied. Section 3 presents our new approach relying on sensitivity analysis, and our novel probabilistic upper bound on the JSR in the quadratic case. In Section 4, we leverage this approach in order to obtain our bound relying on SOS functions. Finally, we illustrate the performance of our bounds in Section 5 with a numerical experiment.
\\
The notation used in this paper is as follows.
For a square matrix $P$, $P \succ $ (resp. $\succeq$) $0$ states that $P$ is positive definite (resp. semi-definite).
We denote the set of real symmetric matrices
of size $n$ by $\mathcal{S}^n$, and the set of positive definite matrices by $\mathcal{S}^n_{++}$.  
For $P \in \mathcal{S}^n_{++}$, we define the \emph{ellipsoidal vector norm} associated to $P$ as: $ \ x \longrightarrow ||x||_P:=\sqrt{x^TPx}$.  $||x||:=\sqrt{x^Tx}$ denotes the classical $\ell_2-$norm of $x \in \mathbb{R}^n$.
\\
Our notation concerning measure-theoretic concepts is taken from \citet{Led:91}. 
$\mathbb{S}=\{x \in \mathbb{R}^n: ||x||=1\}$ denotes the unit sphere centered at the origin in $\mathbb{R}^n$. 
$\sigma^{n-1}$ denotes the \emph{uniform spherical measure}, satisfying $\sigma^{n-1}(\mathbb{S})=1$.
$\mu_M$ denotes the uniform measure on $M=\{1, 2,  . . . , m\}$ for $m \in \mathbb{N}_0$, and $\mu_{M^l}=\otimes^l\mu_M$ denotes the uniform product measure on $M^l=\otimes^l M$, for $l \in \mathbb{N}_0$. Finally, $\mu_l = \sigma^{n-1} \otimes \mu_{M^l}$ denotes the uniform measure on $Z_l=\mathbb{S} \times M^l$. 
\section{Preliminaries}
\subsection{Approximation of the JSR using quadratic forms}
A framework to approximate the JSR of switched linear systems and hence to analyze their stability is provided by CQLFs:
\begin{prop}
(\cite{Jun:09}, Prop. 2.8)
\\
Consider a finite set of matrices $\mathcal{M} \in \mathbb{R}^{n \times n}$.
 If there exists $\gamma \geq 0$ and $P \in \mathcal{S}^n_{++}$ such that  \begin{equation}
     \forall A \in \mathcal{M}, A^TPA \preceq \gamma^2 P
     \label{CQF}
 \end{equation} then $\rho(\mathcal{M}) \leq \gamma$.
\label{theo1}
\end{prop}
When $\gamma<1$, $\|\cdot\|_P$ is a CQLF and can serve as stability certificate for System (\ref{system}). In the spirit of \cite{Ken:19} and in order to benefit from the observation of long trajectories, we consider CQFLs for the finite set $\mathcal{M}^l:=\{\prod_{j=1}^lA_{i_j}:\ A_{i_j} \in \mathcal{M} \}$ and use the relation $\rho(\mathcal{M}^l)=\rho(\mathcal{M})^l$ (\cite{Jun:09}, Prop. 2.5). Since the quality of the upper bound in Proposition \ref{theo1} increases when $\gamma$ decreases, we seek for the minimal $\gamma$ such that \eqref{CQF} is satisfied, by considering the following optimization problem:
\begin{align}
\begin{split}
&\underset{P, \gamma \geq 0}{\min \ } \quad \gamma  \\
& \text{ s.t. } \quad (\mathbf{A}x)^TP\mathbf{A}x\leq \gamma^{2l} x^TPx ,\ \forall \mathbf{A} \in \mathcal{M}^l, \forall x \in \mathbb{S} \\
&\quad \quad \quad \  P \succeq I , \quad \quad ||P||\leq C \text{ for a large } C  \in \mathbb{R}_{\geq 0}.
\label{quadra}
\end{split}
\end{align}
We denote by $(\gamma^o,P^o)$ the solution to Problem \eqref{quadra}. While
the optimal cost $\gamma^o$
is unique, there can be several optimal
$P$. When this occurs, we apply a \emph{tie-breaking rule}: $P^o$
denotes the optimal $P$ with the smallest condition number, since the quality of our upper bound will depend on this number.  Problem \eqref{quadra} differs from Program \eqref{CQF} in three ways. First, we impose a bound on the norm of $P$ to ensure the compactness of the set of constraints. Second, $P \succeq I$ replaces the constraint $P \succ 0$, and third, we restrict $x$ to $\mathbb{S}$. These latter restrictions do not incur any conservativeness because of the \emph{homogeneity} of System (\ref{system}): for $\lambda >0$, $A_{\tau(k)}\lambda x_k=\lambda A_{\tau(k)}x_k$, implying for instance that the decrease of a CQLF on an arbitrary set enclosing the origin is sufficient to serve as a stability certificate.
\\
The bound on the JSR obtained by Program \eqref{quadra} can, in the white-box setting, be improved by considering SOS polynomials $p(x)$, that is polynomials admitting a decomposition $p(x)=\sum_i p_i(x)^2$ (\cite{Par:08}). Before we investigate on how this generalizes to the black-box setting, we formally present the problem addressed in this paper as well as the new data-driven approach we propose in the quadratic case. 
\subsection{Problem statement}
Problem (\ref{quadra}) is defined by an infinite number of constraints, while we only sample a finite number $N$ of trajectories of length $l$ $(x_{i,0}, x_{i,1}, . . . , x_{i,l})$ of System (\ref{system}), corresponding to a  finite subset of the constraints of Problem \eqref{quadra}. These trajectories as well as an upper bound $m$ on the number of modes of the actual system are the only available information about the system.
\\
The trajectories are assumed to be generated from $N$ initial states $x_{i,0}$ drawn randomly, uniformly and independently from  $\mathbb{S}$ according to $\sigma^{n-1}$. Then, $N$ sequences of $l$ modes are drawn randomly, uniformly and independently from $M^l$ according to $\mu_M$ and applied to each initial condition. The set of $N$ observations of System \eqref{system} $\{(x_{i,0},x_{i,1},...,x_{i,l}),\ i=1,2,...,N\}$ is therefore associated to a uniform sample of $N$ $(l+1)$-tuples in $Z_l$:
\begin{equation}
    \omega_N := \{(x_{i,0},j_{i,1},...,j_{i,l}), \ i=1,2,...,N\} \subset Z_l.
\end{equation}
We denote \textbf{$\mathbf{A_j}$}$:= A_{j_{l}}A_{j_{l-1}}...A_{j_1}$ and $\textbf{j}:=\{j_1,...,j_l\}$.
\\
For a given data set $\omega_N$, we can define a \emph{sampled optimization problem} associated to ($\ref{quadra}$), that is a similar problem composed of the $N$ constraints associated to the $N$ observed trajectories. We denote it by Opt($\omega_N$):
\begin{align}
\begin{split}
&\underset{P, \gamma \geq 0}{\min \ } \quad \gamma\\
& \text{ s.t.} \ \quad (\mathbf{A_j}x)^TP\mathbf{A_j}x\leq \gamma^{2l} x^TPx ,\quad \forall (x,\textbf{j}) \in \omega_N\\
& \quad \quad \quad \ P \succeq I,  \quad  \quad ||P|| \leq C \text{ for a large } C  \in \mathbb{R}_{\geq 0}.
\label{quadra1}
\end{split}
\end{align}
Its solution is denoted by $(\gamma^*(\omega_N),P^*(\omega_N))$, where the tie-breaking rule is applied for $P^*(\omega_N)$ if needed. Observe that Program \eqref{quadra1} only requires the knowledge of $\mathbf{A_j}x$, available through the observations: it does not require the knowledge of $\bf{j}$, which is not observed in our setting. Once the solution to Problem \eqref{quadra1} is obtained, we aim to infer the solution to Problem \eqref{quadra} within a certain confidence.
\section{A new data-driven approach for the JSR computation}
In this section, we provide a sensitivity analysis of the quasi-convex optimization problem \eqref{quadra}, allowing to approximate its solution. From this, we obtain a novel probabilistic upper bound on the JSR of $\mathcal{M}$. An asset of this approach is the possibility to extend it to the SOS framework and hence to take advantage of the tighter approximation of the JSR allowed by SOS Lyapunov functions.
\subsection{Sensitivity analysis}
First, we derive a key result about the existence and the cardinality of a finite set of points such that the solutions to Problems \eqref{quadra} and \eqref{quadra1} defined on this set  coincide\footnote{The fundamental difference between Lemma \ref{lem1} and Theorem 2 in \cite{Ber:21} lies in the fact that our bound is derived for Problem \eqref{quadra} with infinite number of constraints while the bound of \cite{Ber:21} is derived for the sampled problem \eqref{quadra1} with a finite number of constraints.}.
\begin{lem}
Let $\gamma^o$ be the optimal cost of Problem (\ref{quadra}).
Then, 
there exists a set $\omega_s \subset Z_l$ with $|\omega_s|=\frac{n(n+1)}{2}+1$ such that $\gamma^*(\omega_s)=\gamma^o$, where $\gamma^*(\omega_s)$ is the optimal solution to Opt ($\omega_s$).
\label{lem1}
\end{lem}
\begin{pf}
From a standard continuity argument (see Appendix A), it is equivalent to prove that for any $\epsilon>0$,
there exists a set $\omega_s \subset Z_l$ with $|\omega_s|=\frac{n(n+1)}{2}+1$ such that $\gamma^*(\omega_s)>\gamma^o-\epsilon$. Let $d_1=\frac{n(n+1)}{2}$ be the number of variables of $P \in \mathcal{S}^n$. Problem (\ref{quadra}) is equivalent to:
\begin{align}
\begin{split}
&\underset{P, \gamma \geq 0}{\min \ } \gamma, \    \text{ s.t. } \quad P \in \mathcal{X}_{(x,\mathbf{j},\gamma)}, \ \forall (x,\mathbf{j}) \in Z_l
\label{quadrabis}
\end{split}
\end{align}
where
$\mathcal{X}_{(x,\mathbf{j},\gamma)}=\{P \in \mathcal{S}^n|\ (\mathbf{A_j}x)^TP\mathbf{A_j}x\leq \gamma^{2l} x^TPx \text{, } P \succeq I \text{, }||P|| \leq C\} \subseteq \mathbb{R}^{d_1}$. Clearly, all sets $\mathcal{X}_{(x,\mathbf{j},\gamma)}$ are convex and compact. Suppose by contradiction that there exists an $\epsilon>0$ such that, for any $\omega \in Z_l$ with $|\omega|=\d_1+1$, it holds that $\gamma^*(\omega)\leq \gamma^o-\epsilon$. Let $\gamma'=\gamma^o-\epsilon$. 
Then, for any $\omega \subset Z_l$ with $|\omega|=\d_1+1$, $\cap_{(x,\mathbf{j}) \in \omega} \mathcal{X}_{(x,\mathbf{j},\gamma')}$ is nonempty. Hence, by Helly's theorem (\cite{Dan:63}), 
 $\cap_{(x,\mathbf{j}) \subset Z_l} \mathcal{X}_{(x,\mathbf{j},\gamma')}$ is nonempty. As a result, $\gamma^o\leq \gamma'=\gamma^o-\epsilon$, which leads to a contradiction.
\end{pf}
Based on Lemma \ref{lem1}, we will show how to construct, from a set of trajectories $\omega_N$, another set $\omega'_N$ such that $\gamma^*(\omega'_N) = \gamma^o$ with a certain probability depending on $N$ and on how close the points of both sets are. Firstly, we introduce the notion of a \emph{spherical cap}:
 \begin{defn}
(\cite{Li:11})
\\
The spherical cap on $\mathbb{S}$ of direction $c$ and measure $\epsilon$ is defined as $\mathcal{C}(c,\epsilon):=\{x \in \mathbb{S}: \ c^Tx>||c||\delta(\epsilon)\} $, where $\delta(\epsilon)$ is given by:
\begin{equation}
    \delta(\epsilon)=\sqrt{1-I^{-1}(2\epsilon; \frac{n-1}{2},\frac{1}{2})}
    \label{delta}
\end{equation}
where $I^{-1}(y;a,b)$ is the \emph{inversed regularized incomplete beta function} (\cite{Maj:73}). Its output is $x>0$ such that $I(x;a,b)=y$, where $I$ is given by
    \[ I :
  \begin{cases}
    \mathbb{R}_{>0} \times \mathbb{R}_{>0} \times \mathbb{R}_{>0} \longrightarrow \mathbb{R}_{\geq 0}    \\
    (x,a,b) \longrightarrow I(x;a,b)=\frac{\int_0^x t^{a-1} (1-t)^{b-1} dt}{\int_0^1 t^{a-1} (1-t)^{b-1} dt}.
  \end{cases}
\]
\end{defn}
\begin{prop}
Let $\omega_N:=\{(x_i,\textbf{j}_i),\ i=1,...,N\}$ be a uniform random sample drawn from $Z_l$, with $N \geq \frac{n(n+1)}{2}+1$, and $\gamma^o$ be the optimal cost of Problem (\ref{quadra}).\\
For all $\epsilon \in (0,1]$, with probability at least $$\beta(\epsilon, m, N):=1-(\frac{n(n+1)}{2}+1)(1-\frac{\epsilon}{m^l})^N,$$there exists a set  $\omega'_N=\{(x'_i,\textbf{j}_i),\ i=1,...,N\}$ such that
\begin{itemize}
    \item $\gamma^*(\omega'_N) = \gamma^o$, where $\gamma^*(\omega'_N)$ is the optimal solution to Opt($\omega_N'$).
    \item For $i=1,...,N$: $x'_i \in \mathbb{S}$ and $||x_i-x'_i||\leq \sqrt{2-2\delta(\epsilon)}$, where $\delta(\cdot)$ is given by (\ref{delta}).
\end{itemize}
\label{omegaprime}
\end{prop}
\begin{pf}
By Lemma \ref{lem1}, there exists a set $\omega_s$, with $|\omega_s| = \frac{n(n+1)}{2}+1$, such that $\gamma^*(\omega_s)=\gamma^o$. Given any $(x',\textbf{j'}) \in \omega_s$, let $(x,\textbf{j})$ be drawn randomly and uniformly from $Z_l$. Then, the probability that $x \in \mathcal{C}(x',\epsilon)$ and $\textbf{j}=\textbf{j'}$ is of $\frac{\epsilon}{m^l}$. Hence, $\mu_l\{(x,\mathbf{j})|\ x \notin \mathcal{C}(x',\epsilon) \text{ or }\textbf{j} \neq \textbf{j'} \}=1-\frac{\epsilon}{m^l}$ for all $(x,\textbf{j}) \in \omega_N$, and $\mu_l^N\{\omega_N \in Z_l^N|\ x_i \notin \mathcal{C}(x',\epsilon) \text{ or }\textbf{j}_i \neq \textbf{j'} \ \forall (x_i,\mathbf{j}_i) \in \omega_N\}=(1-\frac{\epsilon}{m^l})^N$. Since these sets, for each $(x',\textbf{j'}) \in \omega_s$, can be disjoint; the measure of their union, that is the probability that there exists no $(x,\textbf{j}) \in \omega_N$ such that $x \in \mathcal{C}(x',\epsilon)$ and $\textbf{j}=\textbf{j'}$ for a least one $(x',\textbf{j'}) \in \omega_s$, is smaller or equal to $(\frac{n(n+1)}{2}+1)(1-\frac{\epsilon}{m^l})^N$. Therefore, the probability that $\forall (x',\textbf{j'}) \in \omega_s$, there exists a point $(x,\textbf{j}) \in \omega_N$ such that $x \in \mathcal{C}(x',\epsilon)$ and $\textbf{j}=\textbf{j'}$ is at least of $ \beta(\epsilon, m, N)=1-(\frac{n(n+1)}{2}+1)(1-\frac{\epsilon}{m^l})^N$.
\\
If $x \in \mathcal{C}(x',\epsilon)$, then $x'^Tx>||x'||\delta(\epsilon)$ where $\delta(\epsilon)$ is given by \eqref{delta}, and $||x- x'||^2=||x||^2+||x'||^2-2x^Tx'\leq 2-2\delta(\epsilon)$.
We define $\omega'_N$ from $\omega_N$ as follows: for the $(x,\textbf{j}) \in \omega_N$ such that $x \in \mathcal{C}(x',\epsilon)$ and $\textbf{j}=\textbf{j'}$ for some $(x',\mathbf{j'}) \in \omega_s$, $(x',\mathbf{j'}) \in \omega'_N$. For all other $(x,\textbf{j}) \in \omega_N$, $(x,\textbf{j}) \in \omega'_N$.
Hence, with probability larger than $\beta(\epsilon, m, N)$, $\omega_N'$ contains $\omega_s$, which implies $\gamma^*(\omega_N')=\gamma^*(\omega_s)=\gamma^o$. Moreover, by construction, for $1=1,...,N$,  $x'_i \in \mathbb{S}$ and $||x_i-x'_i||\leq \sqrt{2-2\delta(\epsilon)}$.
\end{pf}
Given $\beta, \ m, \ N$, we define $\epsilon(\beta, m,N)=m^l\left(1-\sqrt[N]{\frac{2(1-\beta)}{n(n+1)+2}}\right)$.
\subsection{A new bound for the quadratic case}
Using Proposition \ref{omegaprime}, we provide a probabilistic upper bound on $\gamma^o$, based on $(\gamma^*(\omega_N),P^*(\omega_N))$, the only solution we have access to:
\begin{prop}
Let $\omega_N$  be a uniform random sample drawn from $Z_l$, with $N \geq \frac{n(n+1)}{2}+1$; $(\gamma^*(\omega_N),P^*(\omega_N))$ be the optimal solution to Problem (\ref{quadra1}), and $(\gamma^o,P^o)$ be the optimal solution to Problem (\ref{quadra}). For any $\beta \in [0,1)$, let $\epsilon=m^l\left(1-\sqrt[N]{\frac{2(1-\beta)}{n(n+1)+2}}\right)$. Then, with probability at least $\beta$, the following holds:
\begin{equation}
    \gamma^o\leq \sqrt[l]{ \gamma^*(\omega_N)^l+ (\gamma^*(\omega_N)^l+  B(\mathcal{M}^l))\Delta(\epsilon) \kappa(P^*(\omega_N))}
    \label{rela}
\end{equation}
where $\Delta(\epsilon)=\sqrt{2-2\delta(\epsilon)}$, $\delta(\cdot)$ is given by (\ref{delta}), $B(\mathcal{M}^l)=\underset{\mathbf{A} \in \mathcal{M}^l}{\max \ }||\mathbf{A}||$ and $\kappa(P)=\sqrt{\frac{\lambda_{\max}(P)}{\lambda_{\min}(P)}}$.
\label{rel}
\end{prop}
\begin{pf}
Let $\gamma=\gamma^*(\omega_N)$ and $P=P^*(\omega_N)$. By definition,
    \begin{align}
        ||\mathbf{A_j}x||_P\leq \gamma^l ||x||_P\ \forall (x,\textbf{j}) \in \omega_N.
        \label{base}
    \end{align}
    Consider the Cholesky decomposition $P=L^TL$. Then (\cite{Gol:89}):
    \begin{equation}
     ||x||_P=\sqrt{x^TL^TLx}=||Lx|| \text{ and } ||L^{-1}x||_P=||x||.
        \label{norm1}
    \end{equation}
    Since $
 ||A||_P=||LAL^{-1}||$, $||L^{-1}||_P=||L^{-1}||$. Hence,
    \begin{align}
    \begin{split}
 ||L||||x||\underset{\eqref{norm1}}{\leq}||L||||L^{-1}||||x||_P=\kappa(P)||x||_P.
        \label{norm4}
        \end{split}
    \end{align} 
    Finally,
    \begin{align}
        ||x||_P = ||(x+\Delta x)-\Delta x||_P\leq ||x+\Delta x||_P+||\Delta x||_P.
        \label{zhem}
    \end{align}
    So, for any $\Delta x$ such that $||\Delta x||\leq \Delta(\epsilon)$ and $x+\Delta x \in \mathbb{S}$, we obtain that $\forall (x,\textbf{j}) \in \omega_N$
    \begin{align*}
        &||\mathbf{A_j}(x+\Delta x)||_P\leq ||\mathbf{A_j}x||_P+||\mathbf{A_j}\Delta x||_P\\
        &\underset{(\ref{base}), \ (\ref{norm1})}{\leq} \gamma^l ||x||_P + ||L\mathbf{A_j}\Delta x||\\
        &\underset{\eqref{zhem}, \eqref{norm1}}{\leq } \gamma^l ||x+\Delta x||_P+  (\gamma^l +||\mathbf{A_j}||) ||L||\  ||\Delta x||\\
        & \underset{}{\leq} \gamma^l ||x+\Delta x||_P+ (\gamma^l+  ||\mathbf{A_j}||)\Delta(\epsilon)||L||||x+\Delta x||\\
        & \underset{\eqref{norm4}}{\leq} \left(\gamma^l+ (\gamma^l+  B(\mathcal{M}^l))\Delta(\epsilon) \kappa(P)\right)||x+\Delta x||_P.
    \end{align*}
    Let $\Delta(\epsilon)= \sqrt{2-2\delta(\epsilon)}$, $\beta \in [0,1)$, and $\epsilon(\beta, m,N)=m^l\left(1-\sqrt[N]{\frac{2(1-\beta)}{n(n+1)+2}}\right)$. Then, Proposition \ref{omegaprime} guarantees the existence of a set $\omega'_N=\{(x_i+\Delta x_i,\textbf{j}_i),\ i=1,...,N\}$ such that for $i=1,...,N$, $\Delta x_i \leq \Delta(\epsilon)$ and $x_i+\Delta x_i \in \mathbb{S}$, and such that $\gamma^*(\omega'_N) = \gamma^o$ with probability at least $\beta$.
    \\
    Hence, $||\mathbf{A_j}x||_P\leq \left(\gamma^l+ (\gamma^l+  B(\mathcal{M}^l))\Delta(\epsilon) \kappa(P)\right)||x||_P$ $\forall (x,\textbf{j}) \in \omega'_N$. From the definition of $\gamma^*(\omega'_N) $, $(\gamma^*(\omega'_N))^l\leq \left(\gamma^l+ (\gamma^l+  B(\mathcal{M}^l))\Delta(\epsilon) \kappa(P)\right)$, which concludes
    the proof.
\end{pf}
In the bound of Proposition \ref{rel}, the expression $B(\mathcal{M}^l)=\underset{\mathbf{A} \in \mathcal{M}^l}{\max \ }||\mathbf{A}||$ appears. It can be computed as (\cite{Jun:09}, Prop. 2.7): 
\begin{align}
    B(\mathcal{M}^l)
    =\underset{\lambda \geq 0}{\min \ }\lambda, \text{ s.t. } ||\mathbf{A}x||\leq \lambda,\  \forall x \in \mathbb{S}, \forall \mathbf{A} \in \mathcal{M}^l.
    \label{realmax}
\end{align}
Since we have no direct knowledge of $\mathbf{A} \in \mathcal{M}^l$, we seek to obtain a probabilistic bound on $B(\mathcal{M}^l)$ by solving the sampled problem
\begin{align}
    \underset{\lambda \geq 0}{\min \ } \lambda, \text{ s.t. } ||\mathbf{A_j}x||\leq \lambda ,  \ \forall (x,\mathbf{j}) \in \omega_N
    \label{maxA}
\end{align}
which only requires the knowledge of the observations we have access to. Let $\lambda^*(\omega_N)$ denote the solution to Problem \eqref{maxA}. By combining classical results from chance-constrained optimization and the approach proposed in \cite{Ken:19}, we obtain the following bound: 
\begin{prop}
Consider a set of matrices $\mathcal{M}$, and a uniform random sampling $\omega_N \subset Z_l$, where $N>2$. Let $\lambda^*(\omega_N)$ be the optimal solution to Problem (\ref{maxA}). Then, for any given $\beta_1 \in [0,1)$, with probability at least $\beta_1$, we have:
\begin{equation}
    \underset{\mathbf{A} \in \mathcal{M}^l}{\max \ }||\mathbf{A}|| \leq \frac{\lambda^*(\omega_N)}{\sqrt[l]{\delta(\epsilon_1(\beta_1,m,N))}}
    \label{max2}
\end{equation}
    where $\delta(\cdot)$ is given by (\ref{delta}) and $\epsilon_1(\beta_1,m,N)=\frac{m^l}{2}(1-\sqrt[N]{1-\beta_1})$.
\label{max}
\end{prop}
A proof of Proposition \ref{max} is provided in the Appendix.
We now propose our novel bound, by combining Propositions \ref{rel} and \ref{max}.
\begin{thm}
Consider an $n$-dimensional switched linear system as in (\ref{system}) and a uniform random sampling $\omega_N \subset Z_l$, where $N \geq \frac{n(n+1)}{2}+1$. Let $(\gamma^*(\omega_N), P^*(\omega_N))$ be the optimal solution to (\ref{quadra1}), and $\lambda^*(\omega_N)$ be the optimal solution to (\ref{maxA}). For any $\beta \in [0,1)$ and $\beta_1 \in [0,1)$, let $\epsilon=m^l\left(1-\sqrt[N]{\frac{2(1-\beta)}{n(n+1)+2}}\right)$, and 
$\epsilon_1=\frac{m^l}{2}(1-\sqrt[N]{1-\beta_1})$.  Then, with probability at least $\beta+\beta_1-1$, we have:
\begin{align}
    \rho(\mathcal{M})\leq \sqrt[l]{ \gamma^*(\omega_N)^l+ (\gamma^*(\omega_N)^l+  A(\epsilon_1))\Delta(\epsilon) \kappa(P^*(\omega_N))},
    \label{long}
\end{align}
where $\Delta(\epsilon)=\sqrt{2-2\delta(\epsilon)}$, $\delta(\cdot)$ is given by (\ref{delta}) and $A(\epsilon_1)=\frac{\lambda^*(\omega_N)}{\sqrt[l]{\delta(\epsilon_1)}}$.
    \label{main}
\end{thm}

\begin{pf}
By Proposition \ref{rel}, $\mu_l^N\{\omega_N \in Z_l^N|\ \eqref{rela} \text{ holds}\} \geq \beta$, and by Proposition \ref{max}, $\mu_l^N\{\omega_N \in Z_l^N|\ \eqref{max2} \text{ holds}\}\geq \beta_1$. Hence,
\begin{align*}
    &\mu_l^N\{\omega_N \in Z_l^N|\ \text{at least one of } \eqref{rela} \text{ or } \eqref{max2} \text{ does not hold}\}\\&\leq \mu_l^N\{\omega_N \in Z_l^N|\ \eqref{rela} \text{ does not hold}\}\\
    &\quad + \mu_l^N\{\omega_N \in Z_l^N|\  \eqref{max2} \text{ does not hold}\}\\
    &\leq (1-\beta)+(1-\beta_1).
\end{align*} Therefore, $\mu_l^N\{\omega_N \in Z_l^N|\ \eqref{rela} \text{ and } \eqref{max2} \text{ hold}\}\geq 1- (1-\beta)-(1-\beta_1)$.
If $\omega_N$ is such that $\eqref{rela} \text{ and } \eqref{max2} \text{ hold}$, then \eqref{long} holds. Hence, this relation holds with probability at least $\beta+\beta_1-1$, which concludes the proof.
\end{pf}
\begin{rem}
$\delta(\cdot)$ is defined over $(0,\frac{1}{2})$. If $\epsilon$ is larger than $\frac{1}{2}$, $\delta(\epsilon)=0$. Consequently, our upper bound presents an advantage over the bound in \cite{Ken:19}: in the latter, especially when the number of modes $m$ and the length of the traces $l$ increase, a larger number of samples $N$ is required to provide a finite upper bound on the JSR. Indeed, in our bound, if $\epsilon(\beta, m, N) \geq \frac{1}{2}$ so that $\delta(\epsilon(\beta, m, N))=0$, the bound we provide is nevertheless finite. If $\epsilon_1(\beta_1,m,N)\geq \frac{1}{2}$, the only bound we can provide is $+\infty$, but given that $\epsilon_1(\beta_1,m,N)\leq\frac{\epsilon(\beta, N)m^l}{2}\sqrt{\frac{det(P^*(\omega_N))}{\lambda_{\min}(P^*(\omega_N))^n}}$, a smaller number of samples is required to obtain a finite upper bound. For example, for $m=5$, $n=2$ and $l=3$, our bound is finite from $N=375$ samples while
$N=500$ samples are required to obtain a finite upper bound in \cite{Ken:19}.
\end{rem}
\section{Probabilistic SOS-based guarantees}
We will now extend our approach to the SOS framework. First, we recall properties of SOS Lyapunov functions in the white-box setting, and present generalized versions of Problems \eqref{quadra} and \eqref{quadra1}. Then, we derive a probabilistic upper bound on the JSR using sensitivity analysis. 
Before proceeding further,  we recall the notion of \emph{symmetric algebra of a vector space}:
\begin{defn}
(\cite{Jun:09} and \cite{Par:08})
\\
Let $x \in \mathbb{R}^{n}$. Let $D$ denote the number of different monomials of degree $d$:
\begin{align}
    D=\binom{n+d-1}{d}.
    \label{D}
\end{align}
\\
The \emph{$d$-lift} of x, denoted by $x^{[d]}$, is the vector in $\mathbb{R}^D$, indexed by all the possible exponents $\alpha$ of degree $d$:
\begin{align}
    x^{[d]}_{\alpha}=\sqrt{\alpha !}x^{\alpha},
    \label{dlift}
\end{align}
where $\alpha=(\alpha_1,...,\alpha_n)$ with $|\alpha|=d$, and $\alpha!=\binom{d}{\alpha_1,...,\alpha_n}$ is the multinomial coefficient.
The $d$-lift of the matrix $A$ is the matrix $A^{[d]} \in \mathbb{R}^{D\times D}$ associated to the linear map
$A^{[d]}: x^{[d]} \longrightarrow (Ax)^{[d]}$.
\end{defn}
With the lift as defined in \eqref{dlift}, one has:
\begin{equation}
    ||x^{[d]}||=||x||^d.
    \label{xd}
\end{equation}
We are now ready to introduce SOS forms to approximate the JSR of $\mathcal{M}$.
\begin{prop}
(\cite{Par:08}, Thm. 2.3)
\\
A homogeneous multivariate polynomial $p(x)$ of degree $2d$ is a SOS polynomial if and only if for some $P \in \mathcal{S}^D_{++}$, where $D$ is given by \eqref{D},  $$p(x)=(x^{[d]})^TPx^{[d]}.$$
\end{prop}
We can obtain tighter approximations of the JSR using SOS functions due to the following property:
\begin{prop}
(\cite{Jun:09}, Thm. 2.13)
Consider a finite set of matrices $\mathcal{M} \in \mathbb{R}^{n \times n}$.
 If there exists $\gamma \geq 0$ and $P \in \mathcal{S}^D_{++}$ such that $$\forall \mathbf{A} \in \mathcal{M}^l, \ \forall x \in \mathbb{S}, \ ((\mathbf{A}x)^{[d]})^TP(\mathbf{A}x)^{[d]} \leq \gamma^{2dl} (x^{[d]})^TPx^{[d]}$$
    then $\rho(\mathcal{M}) \leq \gamma$.
\label{lbsos}
\end{prop}
The restriction of $x$ to $\mathbb{S}$ is again due to the homogeneity of System (\ref{system}). SOS functions hence become SOS Lyapunov functions when $\gamma <1$. We generalize Problem \eqref{quadra} as follows:
\begin{align}
\begin{split}
&\underset{P, \gamma \geq 0}{\min \ } \gamma\\
& \text{ s.t.} (\mathbf{A}x)^{[d]T}P(\mathbf{A}x)^{[d]} \leq \gamma^{2dl} x^{[d]T}Px^{[d]}, \forall \mathbf{A} \in \mathcal{M}^l, \forall x \in \mathbb{S}\\
&\quad \quad P \succeq I
, \quad \quad \quad ||P||\leq C \text{ for a large } C  \in \mathbb{R}_{\geq 0}.
\label{sosprob}
\end{split}
\end{align}
Its solution is denoted by $(\gamma_{sos}^o,P_{sos}^o)$, with the same tie-breaking rule applied to $P_{sos}^o$ when the latter is not unique. When $d=1$, we recover Problem \eqref{quadra}.
We consider the following associated sampled problem, denoted by Opt$_{sos}(\omega_N)$ and with solution $(\gamma^*_{sos}(\omega_N), P^*_{sos}(\omega_N))$:
\begin{align}
\begin{split}
&\underset{P, \gamma \geq 0}{\min \ } \gamma\\
& \text{ s.t.} ((\mathbf{A_j}x)^{[d]})^TP(\mathbf{A_j}x)^{[d]} \leq \gamma^{2dl} (x^{[d]})^TPx^{[d]}, \ \forall (x, \mathbf{j}) \in \omega_N\\
& \quad \quad P \succeq I
, \quad  \quad  \quad ||P|| \leq C \text{ for a large } C  \in \mathbb{R}_{\geq 0}.
\end{split}
\label{sosprob1}
\end{align}
Applying Proposition \ref{omegaprime} to Problem \eqref{sosprob} (only the dimension of the optimization program changes from $\frac{n(n+1)}{2}$ to $\frac{D(D+1)}{2}$), we obtain a similar covering of the sampling space, with a confidence level $\beta(\epsilon_{sos}, m, N)=1-(\frac{D(D+1)}{2}+1)(1-\frac{\epsilon_{sos}}{m^l})^N$. Combining this with Proposition \ref{max},  we obtain a SOS-based probabilistic upper bound on $\rho(\mathcal{M})$:
\begin{thm}
Consider an $n$-dimensional switched linear system as in (\ref{system}). Let $D=\binom{n+d-1}{d}$ be given by \eqref{D}, and $\omega_N \subset Z_l$ be a uniform random sampling, with $N \geq \frac{D(D+1)}{2}+1$. Let $(\gamma^*(\omega_N), P^*(\omega_N))$ be the optimal solution to (\ref{sosprob1}) and $\lambda^*(\omega_N)$ be the optimal solution to (\ref{maxA}). For any $\beta \in [0,1)$ and $\beta_1 \in [0,1)$, let $\epsilon_{sos}=m^l\left(1-\sqrt[N]{\frac{2(1-\beta)}{D(D+1)+2}}\right)$ and $\epsilon_1=\frac{m^l}{2}(1-\sqrt[N]{1-\beta_1})$. Then, with probability at least $\beta+\beta_1-1$, we have:
\begin{align*}
    &\rho(\mathcal{M}) \leq\\
    &  \sqrt[dl]{\gamma(\omega_N)^{*dl}+ ((\gamma(\omega_N)^{*dl}+ A(\epsilon_1)^d) f(d,\epsilon_{sos})) \kappa(P^*(\omega_N))}
\end{align*}
where $A(\epsilon_1)=\frac{\lambda^*(\omega_N)}{\sqrt[l]{\delta(\epsilon_1)}}$,\\ $f(d,\epsilon_{sos})=\sqrt{D}((1+\Delta(\epsilon_{sos}))^d-1- (1-\frac{1}{\sqrt{D}})\Delta(\epsilon_{sos})^d)$, $\Delta(\epsilon_{sos})=\sqrt{2-2\delta(\epsilon_{sos})}$ and $\delta(\cdot)$ is given by (\ref{delta}). 

    \label{mainSOS}
\end{thm}
When $d=1$, this theorem reduces to Theorem \ref{main} as expected. A proof of Theorem \ref{mainSOS} is provided in the Appendix. Since $\epsilon_1$ does not depend on $d$, the minimum number of samples required
to obtain a finite bound is the same as in the quadratic case.
\section{Experimental results}
\begin{figure}[h]
    \centering
    \includegraphics[width=8.4cm]{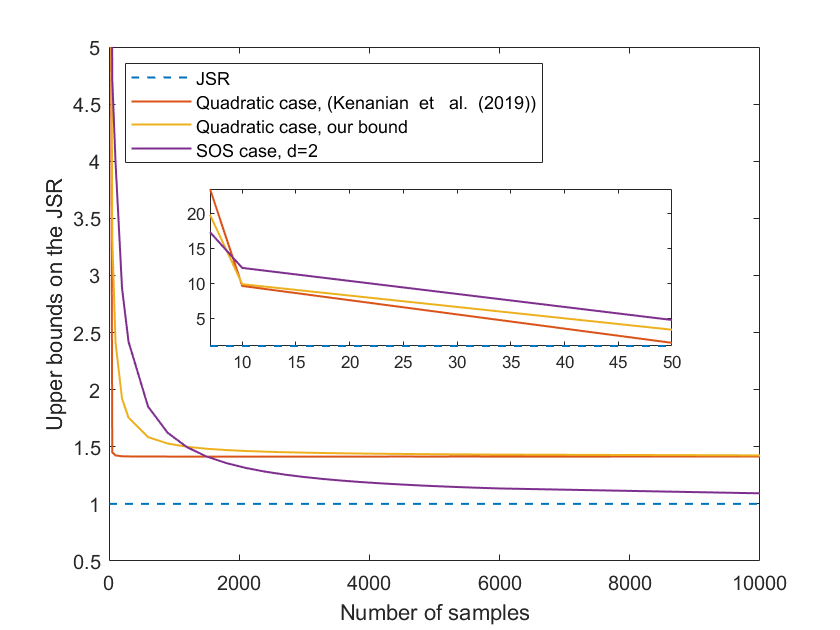}
    \caption{Evolution of the upper bounds with the number of samples, in the quadratic and in the SOS case.}
    \label{SOStest}
\end{figure}

In this section, we conduct an experiment in order to illustrate the improvement of our quadratic bound compared to the one of \cite{Ken:19}, and the performance of our SOS bound compared to the quadratic bounds. We consider a stable system for which no CQLF exists (\cite{Par:08}, Ex. 2.8):  $\mathcal{M}=\{\begin{pmatrix} 1 & 0 \\
1 & 0\end{pmatrix}, \begin{pmatrix} 0 & 1 \\
0 & -1\end{pmatrix}\}$. The true JSR of this set is equal to 1, but the best quadratic bound (as given by \eqref{quadra}) is equal to $\sqrt{2}$, while one can obtain better bounds using SOS functions.
First, using $N$ observations, for N ranging between $10$ and $10000$, we solve Problems (\ref{maxA}),  (\ref{quadra1}) in the quadratic case, and (\ref{sosprob1}) in the SOS case, for a SOS function of degree $2d=4$. Once $\gamma^*(\omega_N)$ is obtained, we apply the tie-breaking rule in order to obtain a feasible $P$ of minimal condition number. Then, we compute the bounds of Theorems \ref{main} and \ref{mainSOS} as well as the bound of \cite{Ken:19}. $\beta$ and $\beta_1$ are fixed to $0.95$, for a global confidence level of $\beta+\beta_1-1=0.9$, and the length of the trace is set to $1$. The average of the bounds is computed over 10 runs for each value of $N$.
\\
Fig. \ref{SOStest} presents the evolution of the different upper bounds with the number of samples. We observe that for a large enough number of samples, the SOS bound outperforms the quadratic bound in the black-box setting just like it does in the white-box setting. We also observe that while the bound of \cite{Ken:19} is tighter than our quadratic bound for a sufficiently high number of samples, our bound is better for a small number of samples. The fact that our bound is better for a small number of samples makes it even more advantageous when longer traces are used. Indeed, for larger traces, because the probability space to be sampled is larger, the performance of the bounds decreases, but faster for the bound of \cite{Ken:19} than for our bound.
\\
\section{Conclusion}
We have introduced a new data-driven approach providing probabilistic upper bounds on the JSR of an unknown system, based on a finite number of observations of trajectories. Our approach relies on sensitivity analysis of a quasi-convex optimization problem, related to the construction of approximate Lyapunov functions. Based on this sensitivity analysis, we firstly derived a new bound using common quadratic Lyapunov functions, and then generalized our method to SOS Lyapunov functions. This approach is well known to improve the quality of approximation of the JSR in the white-box setting. We demonstrated on a numerical experiment that it indeed improves the JSR approximation in the black-box setting too.
\\
For future work, we plan to investigate on the effect of adaptive sampling on the performance of our bounds: we believe that our analysis can be improved by sampling points a posteriori; and to compare our results with the
ones that can be obtained by first identifying the model. We also leave for further work the generalisation of data-driven stability analysis of switched linear systems to other Lyapunov functions, for instance to the path-complete Lyapunov framework.
\bibliography{ifacconf}             

\appendix
\section{Continuity argument}
We prove the equivalence between the two following statements:
\begin{enumerate}
    \item[(B1)] There exists a set $\omega \subset Z_l$ with $|\omega|=d_1+1$ such that $\gamma^*(\omega)=\gamma^o$.
    \item[(B2)] For any $\epsilon>0$,
there exists a set $\omega \subset Z_l$ with $|\omega|=d_1+1$ such that $\gamma^*(\omega)>\gamma^o-\epsilon$.
\end{enumerate}
Clearly, if (B1) holds then (B2) holds. Suppose that (B2) holds. 
Then, we can construct a sequence $\{(\omega_k,\epsilon_k)\}$ where $\omega_k \subset Z_l$, $|\omega_k|=\frac{n(n+1)}{2}+1$ and $\epsilon_{k+1}\leq \epsilon_k$ such that $\epsilon_k \longrightarrow 0$ and $\gamma^*(\omega_k)>\gamma^o-\epsilon_k$. Taking the limit on both sides yields $\underset{k \longrightarrow \infty}{\lim \ }\gamma^*(\omega_{k})\geq \gamma^o$, implying $\underset{k \longrightarrow \infty}{\lim \ }\gamma^*(\omega_{k})= \gamma^o$ since $\forall \omega_k \in Z_l$, $\gamma^*(\omega_{k})\leq \gamma^o$. 
\\
By compactness of the set $Z_l$, $\{\omega_{k}\}$ (or a subsequence of $\{\omega_{k}\}$) converges to a limit $\bar \omega=\underset{k \longrightarrow \infty}{\lim \ }\omega_{k}$. We now show that that $\gamma^*(\cdot)$ is continuous, which implies that
$\gamma^*(\bar \omega)=\underset{k \longrightarrow \infty}{\lim \ }\gamma^*(\omega_{k})=\gamma^o$ so that (S1) holds.
We prove that for a given $\bar \epsilon >0$ defining $\bar \delta:=\frac{\bar \epsilon^l}{2\kappa(P)\underset{\mathbf{A} \in \mathcal{M}^l}{\max \ }||\mathbf{A}||}$, if $\Delta \omega \in  (\bar \delta \mathbb{B}\times 0) \times ... \times (\bar \delta \mathbb{B}\times 0)$, then $|\gamma^*(\bar \omega+\Delta \omega)-\gamma^*(\bar \omega)|\leq \bar \epsilon$.
\\
Indeed, 
\begin{enumerate}
    \item Consider $x$, $\mathbf{j}$ and $\gamma^*$ such that  $||\Delta x||\leq \bar \delta$ and
$    ||\mathbf{A_j}x||_P\leq \gamma^{*l} ||x||_P$. Then, by a reasoning similar to the derivation of Proposition \ref{rel} and since by definition of Problem \eqref{quadra1}, $\gamma^{*l} \leq \underset{\mathbf{A} \in \mathcal{M}^l}{\max \ }||\mathbf{A}||$, we have: \begin{align*}
    &||\mathbf{A_j}(x+\Delta x)||_P\\
    &\leq \left(\gamma^{*l}+ 2\underset{\mathbf{A} \in \mathcal{M}^l}{\max \ }||\mathbf{A}||\bar \delta \kappa\right)||x+\Delta x||_P
    \\
    &\leq(\gamma^{*l}+ \bar \epsilon^l)||x+\Delta x||_P
\end{align*} 
by definition of $\bar \delta$.
\item Conversely, consider $x$, $\mathbf{j}$ and $\gamma^*$ such that  $||\Delta x||\leq \bar \delta$ and
$    ||\mathbf{A_j}(x+\Delta x)||_P\leq \gamma^{*l} ||x+\Delta x||_P$. Then,
\begin{align*}
    ||\mathbf{A_j}x||_P&\leq ||\mathbf{A_j}(x+\Delta x)||_P+||\mathbf{A_j}\Delta x||_P \\
    &\leq \left(\gamma^{*l}+ (\gamma^{*l}+\underset{\mathbf{A} \in \mathcal{M}^l}{\max \ }||\mathbf{A}||)\bar \delta \kappa(P)\right)||x||_P
    \\&
    \leq(\gamma^{*l}+ \bar \epsilon^l)||x||_P
\end{align*} 
\end{enumerate}  
 Hence, since $|\sqrt[l]{x}-\sqrt[l]{y}|\leq \sqrt[l]{|x-y|}$, the result holds.
\\
Finally, from the convergence of $\{\omega_k\}$, for any $0<\epsilon<1$, there exists $l$ such that the components in $\omega_k$ corresponding to the modes do not change
(or have reached the limit) for all $k>l$, which concludes the proof.
\section{Proof of Proposition \ref{max}}
We define an essential set of a sampled problem as a subset, of minimal cardinality, of constraints whose removal implies a decrease in the cost of the problem. Problem \eqref{maxA} is convex and admits a unique solution. Moreover, there exists a sample $(x,\mathbf{j}) \in \omega_N$ such that $\lambda^*(\omega_N)=||\mathbf{A_j}x||$, that is $(x, \mathbf{j})=\underset{(x,\mathbf{j}) \in \omega_N}{\text{argmax }} ||\mathbf{A_j}x||$, so that one is a bound on the cardinality of its essential sets. Let $V(\omega_N)$ be the set $\{(x,\mathbf{j}) \in Z_l:\ ||\mathbf{A}_jx||>\lambda^*(\omega_N) \}$ of violated constraints. We seek to bound the probability to sample $\omega_N$  such that the measure of this set is larger than $\epsilon_1$.
We consider two cases, depending on whether or not the essential set of Problem \eqref{maxA} is unique.
\\
Suppose the essential set is unique with probability one. Then, a classical result from chance-constrained optimization (\cite[Thm. 3.3]{Cal:10}) ensures that for all $\epsilon_1 \in (0,1]$, we have 
$    \mu^N_l\{\omega_N \in Z_l^N:\ \mu_l(V(\omega_N))\leq \epsilon_1\} \geq \beta_1(\epsilon_1, N)$, 
where  $\beta_1(\epsilon_1, N)=1-(1-\epsilon_1)^{N}$. 
\\
Conversely, suppose that the problem admits distinct essential sets with probability larger than 0, which implies the existence of at least one matrix $\mathbf{A} \in \mathcal{M}^l$ such that $\mathbf{A}^T\mathbf{A}=\gamma I_n$ for some $\gamma \in \mathbb{R}_{\geq 0}$ (adapted from \cite{Ber:21}).  Let $\lambda^o=||\mathbf{A}^*||$ be the solution to Problem \eqref{realmax}. We consider the two following cases:
\\
(a) $\mathbf{A}^{*T}\mathbf{A}^*=\lambda^oI_n$. Then, $\mu_l(V(\omega_N))$ can only be positive if $\omega_N$ contains no $(x,\mathbf{j})$ such that $\mathbf{A}^*=\mathbf{A_j}$, which happens with probability $(1-\frac{1}{m})^N$. Hence, $\forall \epsilon_1 \in (0,\frac{1}{m}]$,
    $$
         \mu^N_l\{\omega_N \in Z_l^N:\ \mu_l(V(\omega_N))\geq \epsilon_1\} \leq(1-\frac{1}{m})^N \leq (1-\epsilon_1)^N.
    $$
(b) $\mathbf{A}^{*T}\mathbf{A}^*\neq \lambda^oI_n$.
    Let $V'(\lambda)=\mu_l(\{(x,\mathbf{j}) \in Z_l:\ ||\mathbf{A}_jx||>\lambda\})$. Then, since $1-V'(\lambda)$ is the cumulative distribution function of $||\mathbf{A_j}x||$, $V'(\lambda)$ is right continuous and non increasing. Hence, $\forall \epsilon_1 \in (0,1)$, there exists an $\lambda(\epsilon_1) \in \mathbb{R}_{\geq 0}$ satisfying $V'(\lambda(\epsilon_1))= \epsilon_1$, where $\lambda(\epsilon_1)$ is a decreasing function of $\epsilon_1$. In addition,  $\lambda^*(\omega_N)\leq \lambda(\epsilon_1)$ implies that $\omega_N$ contains no $(x,\mathbf{j}) \in \{(x,\mathbf{j}) \in Z_l:\ ||\mathbf{A}_jx||>\lambda(\epsilon_1) \}$, which occurs with probability $(1-\epsilon_1)^N$. Hence, \begin{align*}
        &\quad \mu^N_l\{\omega_N \in Z_l^N:\ \mu_l(V(\omega_N))\geq \epsilon_1\}\\&=\mu^N_l\{\omega_N \in Z_l^N:\ \lambda^*(\omega_N)\leq \lambda(\epsilon_1)\}\leq(1-\epsilon_1)^N.
    \end{align*}
In all cases, $    \mu^N_l\{\omega_N \in Z_l^N:\ \mu_l(V(\omega_N))\leq \epsilon_1\} \geq \beta_1(\epsilon_1, N)$, where  $\beta_1(\epsilon_1, N)=1-(1-\epsilon_1)^{N}$.
This result can replace Theorem 4 in \cite{Ken:19}. Then, the proof follows the same lines as the derivation of Theorem 6 in \cite{Ken:19}, so that details are omitted here.
\section{Proof of Theorem \ref{mainSOS}}
Throughout this section, let $D$ be given by \eqref{D}.
Before proceeding to the proof itself, some additional definitions and properties are needed. 
\begin{prop}
Given $x, \ y \in \mathbb{R}^n$, the \emph{$k$-th Kronecker power} of $x$ is defined as: $x^{\otimes k}=x\otimes x^{\otimes(k-1)}$ and $x^{\otimes 1}=x$, where $\otimes$ denotes the Kronecker product. Then, $\xd=C_d x^{\otimes d}$, where $C_d \in \mathbb{R}^{D \times n^d}$ is a matrix of norm $||C_d||\leq \sqrt{D}$ and $\xd$ is the d-lift of $x$.
\\
Moreover, the following relations hold: \begin{align}
    ||\xk||=||x||^k    \label{kro1}
\end{align}
    \begin{align}
        (x+y)^{\otimes d}&=(x+y) \otimes ... \otimes (x+y) = \sum_{k=0}^d K_{xy}(d,k) \nonumber\\
        &= \xk+y^{\otimes d}+\sum_{k=1}^{d-1} K_{xy}(d,k)
        \label{kro}
    \end{align}
    where $K_{xy}(d,k)$ is the sum of all possible sequences of $d$ Kronecker products composed of $d-k$ copies of $x$ and $k$ copies of $y$. For $k=0,...,d$, $K_{xy}(d,k)$ is composed of $\binom{d}{k}$ terms of norm $||x^{\otimes (d-k)}\otimes y^{\otimes k}||=||x||^{(d-k)}|| y||^k.$

\label{krokro}
\end{prop}
\textbf{Proof of Proposition 7:} The two last relations are well-known. In addition, $x^{\otimes d}$ is composed of $\alpha!$ copies of each monomial of degree $d$ $x^{\alpha}$. Hence, if the only non-zero entries of $\alpha$-th row of $C_d$ are $\frac{1}{\sqrt{\alpha!}}$ and lie at columns associated to monomials $x^{\alpha}$, $C\xk=x^d$. Moreover, $||C_d||^2\leq ||C_d||_F^2=\sum_{\alpha} \alpha! \frac{1}{\alpha!}=D$.
\\

We are now ready to prove Theorem \ref{mainSOS}.
\begin{pf}
Let $\gamma=\gamma_{sos}^*(\omega_N)$ and $P=P_{sos}^*(\omega_N)$. By definition,
    \begin{align}
        ||\Axd||_P\leq \gamma^{dl} ||\xd||_P,\ \forall (x,\textbf{j}) \in \omega_N.
        \label{ba}
    \end{align}
    Let $\psi=||\dxd||_P$. For any $\Delta x$ such that $||\Delta x||\leq \Delta(\epsilon)$ and $x+\Delta x \in \mathbb{S}$, we obtain that $\forall (x,\textbf{j}) \in \omega_N$
    \begin{align*}
        &||\Axdx||_P\underset{Prop. \ref{krokro}}{=}||\Ad C_d (x+\Delta x)^{\otimes d}||_P\\
        &\underset{(\ref{kro})}{\leq}||\Ad C_d \xk ||_P+ ||\Ad C_d (\Delta x)^{\otimes d}||_P\\
        &+|| \sum_{k=1}^{d-1} \Ad C_d K_{x\Delta x}(d,k) ||_P\\
        &\underset{(\ref{ba})}{\leq }
        \gamma^{dl}||\xd||_P + ||(\mathbf{A_j}\Delta x)^{[d]}||_P\\
        &+\sum_{k=1}^{d-1} || \Ad C_d K_{x\Delta x}(d,k) ||_P
        \\
        & \underset{(\ref{norm1})}{\leq} \gamma^{dl}||\dxd-\Delta x^{[d]}-\sum_{k=1}^{d-1} C_d K_{x\Delta x}(d,k) ||_P\\& +||\mathbf{A_j}\Delta x||^d||L||+\sum_{k=1}^{d-1} ||\Ad|| ||C_d|| ||K_{x\Delta x}(d,k) || ||L||
        \\
        & \underset{(\ref{xd})}{\leq} \gamma^{dl}\psi+ (\gamma^{dl}+ ||\mathbf{A_j}||^d)||\Delta x||^d||L||\\
        &+\sum_{k=1}^{d-1} (\gamma^{dl}+||\mathbf{A_j}||^d) \sqrt{D} \binom{d}{k} ||x||^{d-k}||\Delta x||^k ||L||\\
        &\underset{\eqref{norm4}}{\leq} \gamma^{dl}\psi+ (\gamma^{dl}+B(\mathcal{M}^l)^d)\Delta(\epsilon)^d \kappa(P) \psi \\
        &+\left((\gamma^{dl}+B(\mathcal{M}^l)^d) \sqrt{D}\sum_{k=1}^{d-1} \binom{d}{k} \Delta(\epsilon)^k \right) \kappa(P) \psi
        \\
        & \underset{(\ref{norm1}) }{\leq}\left(\gamma^{dl}+(\gamma^{dl}+ B(\mathcal{M}^l)^d)\Delta(\epsilon)^d\kappa(P)\right)\psi+
        \\&((\gamma^{dl}+B(\mathcal{M}^l)^d) \sqrt{D}(\sum_{k=0}^{d} \binom{d}{k} \Delta(\epsilon)^k-1-\Delta(\epsilon)^d) )\kappa(P)\psi
        \end{align*}
        \begin{align*}
        & \underset{(\ref{norm4})}{=}\left(\gamma^{dl}+(\gamma^{dl}+ B(\mathcal{M}^l)^d)\Delta(\epsilon)^d\kappa(P)\right) \psi+\\
        &\left((\gamma^{dl}+B(\mathcal{M}^l)^d) \sqrt{D}\left((1+\Delta(\epsilon))^d-1-\Delta(\epsilon)^d\right)\right) \kappa(P) \psi
    \end{align*}
    where $B(\mathcal{M}^l)=\underset{\mathbf{A}\in \mathcal{M}^l}{\max}||\mathbf{A}||$. Let $\beta \in [0,1)$, $\Delta(\epsilon)=\sqrt{2-2\delta(\epsilon_{sos})}$ and $\epsilon_{sos}(\beta, m,N)=m^l\left(1-\sqrt[N]{\frac{2(1-\beta)}{D(D+1)+2}}\right)$. The extended version of Proposition \ref{omegaprime} guarantees the existence of a set $\omega'_N=\{(x_i+\Delta x_i,\textbf{j}_i),\ i=1,...,N\}$ such that for $i=1,...,N$, $\Delta x_i \leq \Delta(\epsilon)$ and $x_i+\Delta x_i \in \mathbb{S}$; and such that $\gamma_{sos}^*(\omega'_N) = \gamma_{sos}^o$ with probability at least $\beta$.
    \\
    Hence, $||\Axd||_P\leq \big(\gamma^{dl}+
        ((\gamma^{dl}+ B(\mathcal{M}^l)^d) \sqrt{D}$\\ $\left((1+\Delta(\epsilon))^d-1- (1-\frac{1}{\sqrt{D}})\Delta(\epsilon)^d\right) ) \kappa(P)\big) ||\xd||_P$,\\ $\forall (x,\textbf{j}) \in \omega'_N$.
    Following the same reasoning as in Theorem \ref{main} concludes the proof.
\end{pf}
\end{document}